\documentclass[11pt]{amsart}

\newtheorem{theorem}{Theorem}[section]
\newtheorem{lem}[theorem]{Lemma}

\newtheorem{conj}[theorem]{Conjecture}
\theoremstyle{definition}
\newtheorem{definition}[theorem]{Definition}
\newtheorem{example}[theorem]{Example}
\newtheorem{observation}[theorem]{Observation}

\theoremstyle{remark}

\numberwithin{equation}{section}

\newcommand{\ra}{\rightarrow}

\newcommand{\ddt}{\frac{\text{d}}{\text{dt}}}

\newcommand{\lb}{\langle}
\newcommand{\rb}{\rangle}

\newcommand{\mg}{\mathfrak{g}}
\newcommand{\mh}{\mathfrak{h}}

\renewcommand{\mp}{\mathfrak{p}}

\newcommand{\R}{\mathbb{R}}

\newcommand{\ii}{\mathbf{i}}
\newcommand{\jj}{\mathbf{j}}
\newcommand{\kk}{\mathbf{k}}

\begin{document}

\newcommand{\spacing}[1]{\renewcommand{\baselinestretch}{#1}\large\normalsize}
\spacing{1.14}

\title{Invariant metrics on $SO(4)$ with nonnegative curvature}

\author {Kristopher Tapp$^\ast$}

\address{Department of Mathematics\\ Williams College\\
Williamstown, MA 01267}
\email{ktapp@williams.edu}
\thanks{$^\ast$Supported in part by NSF grant DMS--0303326.}



\begin{abstract}
We derive a curvature-variation formula for a path of left-invariant metrics on a compact Lie group, beginning at a bi-invariant metric.  We prove rigidity theorems for paths which remain nonnegatively curved, and we make progress towards a classification of the left-invariant metrics with nonnegative curvature on $SO(4)$.
\end{abstract}

\maketitle

\section{Introduction}
The starting point for constructing all known examples of compact manifolds with positive (or even quasi-positive) curvature is the fact that bi-invariant metrics on compact Lie groups are nonnegatively curved.  In order to generalize this fundamental starting point, we address the question: given a compact Lie group $G$, classify the left-invariant metrics on $G$ which have nonnegative curvature.  The first two cases, $G=SO(3)$ and $G=U(2)$, were completely solved in~\cite{SMALL}.  For higher dimensional groups, new examples could potentially, via familiar quotient constructions, lead to new examples of quasi-positively curved spaces.  

To explore the remaining cases, we first formulate an infinitesimal version of the problem.  That is, we consider a path $h(t)$ of left-invariant metrics on $G$, with $h(0)$ bi-invariant.  The main work of this paper is towards classifying the possibilities for $h'(0)$ such that the path appears (up to derivative information at $t=0$) to remain nonnegatively curved.  We derive restrictions on $h'(0)$ for general $G$, and some specific to $G=SO(4)$.  These infinitesimal restrictions will be translated into global conclusions (about the set of all left-invariant metrics on $G$ with nonnegative curvature) via power-series techniques in the follow-up paper~\cite{SMALL06}.

The author is pleased to thank Craig Sutton, Emily Proctor, Zachary Madden and Nela Vukmirovic for their help with this project.
\section{Curvature Variation formulas}\label{intro}
Let $G$ be a compact Lie group.  Let $h_0$ denote a bi-invariant metric on $G$.  Let $h$ be an arbitrary left-invariant metric on $G$. 
The value of $h$ at the identity $e\in G$ is determined in terms of $h_0$ by some self-adjoint $\Phi:\mg\ra\mg$ defined so that for all $A,B\in\mg$,
$$h(A,B)=h_0(\Phi(A),B).$$
Denote $\lb A,B\rb = h_0(A,B)$, $|A|^2=h_0(A,A)$ and $|A|_h^2=h(A,A)$ for $A,B\in\mg$.

The unnormalized sectional curvature of vectors $Z_1,Z_2\in\mg$ (the Lie algebra of $G$) with respect to $h$ is given by P\"uttmann's formula from~\cite{put}:
\begin{gather}\label{E:put}
k_h(Z_1,Z_2) = (1/2)\lb[\Phi(Z_1),Z_2] + [Z_1,\Phi(Z_2)],[Z_1,Z_2]\rb - (3/4)|[Z_1,Z_2]|^2_h \\
  +\lb B(Z_1,Z_2),\Phi^{-1}(B(Z_1,Z_2))\rb - \lb B(Z_1,Z_1),\Phi^{-1}(B(Z_2,Z_2))\rb,\notag
\end{gather}
where $B(Z_1,Z_2) = (1/2)([Z_1,\Phi(Z_2)]+[Z_2,\Phi(Z_1)])$.

Next, let $h_t$ be a family of left-invariant metrics on $G$, with $h_0$ bi-invariant. Let $\Phi_t:\mg\ra\mg$ be the endomorphism determining $h_t$ in terms of $h_0$, so that $h_t(A,B)=h_0(\Phi_t A, B)$, and let $\Psi:\mg\ra\mg$ denote the derivative of $\Phi_t$ at $t=0$.

With respect to $h_0$, the zero curvature planes are precisely those spanned by commuting pairs of vectors, so we need only study the curvature variation of those planes.
\begin{lem}\label{pp}
Let $X,Y\in\mg$ with $[X,Y]=0$.  Let $k(t)$ denote the unnormalized sectional curvature of $X$ and $Y$ with respect to $h_t$.  Then $k(0)=0$, $k'(0)=0$ and 
$$2k''(0) = |[X,\Psi Y] + [\Psi X,Y]|^2 = \left|\frac{D}{dt}\Big|_{t=0}[\Phi_t X,\Phi_t Y]\right|^2.$$
\end{lem}

\begin{proof}
By P\"uttmann's formula (Equation~\ref{E:put}),
$$
k(t) = \lb B_t(X,Y),\Phi_t^{-1}B_t(X,Y)\rb - \lb B_t(X,X),\Phi_t^{-1}B_t(Y,Y)\rb,
$$
where $B_t(A,B) = (1/2)([A,\Phi_tB]+[B,\Phi_tA])$.  Using primes to denote derivatives with respect to $t$, we have
\begin{eqnarray*}
k'(t) & = & \lb B'_t(X,Y),\Phi_t^{-1} B_t(X,Y)\rb + \lb B_t(X,Y),(\Phi_t^{-1})'B_t(X,Y)+\Phi_t^{-1} B_t'(X,Y)\rb\\
      &  & -\lb B'_t(X,X),\Phi_t^{-1} B_t(Y,Y)\rb - \lb B_t(X,X),(\Phi_t^{-1})'B_t(Y,Y)+\Phi_t^{-1} B_t'(Y,Y)\rb
\end{eqnarray*}
Since $B_0(X,Y)=B_0(X,X)=B_0(Y,Y)=0$, we see that $k'(0)=0$.  Taking one more derivative yields:
\begin{eqnarray*}
k''(0) & = & 2|B_0'(X,Y)|^2-2\lb B_0'(X,X),B_0'(Y,Y)\rb\\
       & = & (1/2)|[X,\Psi Y]+[Y,\Psi X]|^2 - 2\lb[X,\Psi X],[Y,\Psi Y]\rb\\
       & = & (1/2)|[X,\Psi Y]|^2 + (1/2)|[Y,\Psi X]|^2 + \lb[X,\Psi Y],[Y,\Psi X]\rb
                 - 2\lb[X,\Psi X],[Y,\Psi Y]\rb
\end{eqnarray*}
The third and fourth terms above combine by the Bianchi identity:
\begin{eqnarray*}
\lb[X,\Psi Y],[Y,\Psi X]\rb & = & -\lb[X,[Y,\Psi X]],\Psi Y\rb 
                                 =\lb[Y,[\Psi X,X]],\Psi Y\rb \\
                            & = & -\lb[Y,\Psi Y],[\Psi X,X]\rb 
                                 = \lb[X,\Psi X],[Y,\Psi Y].\rb
\end{eqnarray*}
Thus,
$$k''(0) = (1/2)|[X,\Psi Y]|^2 + (1/2)|[Y,\Psi X]|^2 - \lb[X,\Psi Y],[Y,\Psi X]\rb
         = (1/2)|[X,\Psi Y]-[Y,\Psi X]|^2.$$
\end{proof}

It is interesting that $k''(0)\geq 0$; however, this property limits the lemma's usefulness for deciding which variations remain nonnegatively curved.  This problem is remedied by the next lemma, in which the plane under consideration twists along with the changing metrics.  To avoid dependence on higher order derivatives of $\Phi_t$, and for other reasons which will be clear shortly, we restrict $\Phi_t$ to be a certain type of ``straight path''.  More precisely, we call the variation $\Phi_t$ \underline{inverse-linear} if the inverse-path $\Phi_t^{-1}$ is linear, meaning $\Phi_t=(I-t\Psi)^{-1}$.
\begin{lem}\label{kprime}
Assume that $\Phi_t$ is inverse-linear.  Let $X,Y\in\mg$ with $[X,Y]=0$.  Let $\kappa(t)$ denote the unnormalized sectional curvature of $\Phi^{-1}_t X$ and $\Phi^{-1}_t Y$ with respect to $h_t$.  Then $\kappa(0)=0$, $\kappa'(0)=0$, $\kappa''(0)=0$ and 
\begin{eqnarray*}
(1/6)\kappa'''(0) & = & \lb[X,\Psi Y]+[\Psi X,Y],[\Psi X,\Psi Y]\rb
                   +\lb[\Psi X,X],\Psi[\Psi Y,Y]\rb \\
             &   &   -\lb[X,\Psi Y],\Psi[X,\Psi Y]\rb
                   -\lb[X,\Psi Y],\Psi[\Psi X,Y]\rb
                   -\lb[\Psi X,Y],\Psi [\Psi X,Y]\rb.
\end{eqnarray*}
\end{lem}
We omit the proof, which is longer but also more straightforward than the proof of the previous lemma.
\begin{example}\label{EE} Suppose $H\subset G$ is a Lie subgroup with Lie algebra $\mh\subset\mg$.  For $A\in\mg$, let $A^\mh$ and $A^{\mp}$ denote the projections of $A$ onto and orthogonal to $\mh$ with respect to $h_0$.  The variation $\Phi_t(A)=\frac{1}{1+t}A^{\mh} + A^{\mp}$ in inverse-linear and has nonnegative curvature for $t\in[0,1)$.  In this variation, vectors tangent to $\mh$ are gradually shrunk.  The $t=0$ derivative is $\Psi A=-A^{\mh}$.  Lemma~\ref{kprime} yields:$$(1/6)\kappa'''(0)=|[X^\mh,Y^\mh]|^2.$$
This agrees with Eschenburg's formula from~\cite{Eschenburg}, which says the following: With respect to the metric $h_t$ (for fixed $t\in(0,1)$), the plane spanned by $\{A,B\}$ has zero curvature if and only if $[A^\mh,B^\mh]=0$ and $[\Phi_t(A),\Phi_t(B)]=0$.  In other words, the plane spanned by $\Phi^{-1}_t(X)$ and $\Phi^{-1}_t(Y)$ has zero-curvature if and only if $[X,Y]=0$ and $[X^\mh,Y^\mh]=0$.  This compatibility with Eschenburg's conclusion indicates that Lemma~\ref{kprime} uses a more natural setup than Lemma~\ref{pp}.
\end{example}

\begin{example} For $H\subset G$ as before, consider the variation $\Phi_t(A)=\frac{1}{1-t}A^{\mh} + A^{\mp}$, which gradually enlarges $\mh$.  For this variation, $\kappa'''(0)=-6|[X^\mh,Y^\mh]|^2.$  So this variation immediately picks up some negative curvature unless $[X^\mh,Y^\mh]=0$ whenever $[X,Y]=0$.  If $\mh$ is abelian, then $\kappa'''(0)=0$ for all commuting $X,Y$.  In other words, the lemma suggests that enlarging an abelian subalgebra might preserve nonnegative curvature.  Indeed, it is proven in~\cite{GZ} that enlarging an abelian subalgebra as far as $4/3$ always preserves nonnegative curvature.
\end{example}

The inverse-linear parametrization used in Example~\ref{EE} looks natural when re-described as a family of submersions metrics: $(G,h_t)=((G,h_0)\times(H,(1/t)g_0))/H$.  Further evidence in this paper and in~\cite{SMALL06} will showcase the inverse-linear paths as natural, and will suggest:

\begin{conj}[Inverse-Linear Conjecture]\label{ILC} If $h$ is a left-invariant metric with nonnegative curvature on a compact Lie group $G$, then the unique inverse-linear path from any fixed bi-invariant metric $h_0$ to $h_1=h$ is through nonnegatively curved metrics.
\end{conj}
In particular, the conjecture would imply that the nonnegatively curved left-invariant metrics on $G$ form a path-connected space.

\begin{definition}
Let $\Psi:\mg\ra\mg$ be an $h_0$-self-adjoint map.  Let $\Phi_t=(I-t\Psi)^{-1}$ be the corresponding inverse-linear variation.  We refer to $\Psi$ (or to $\Phi_t$) as \underline{infinitesimally nonnegative} if for each $X,Y\in\mg$ with $[X,Y]=0$, there exists $\epsilon>0$ such that $\kappa(t)\geq 0$ for all $t\in[0,\epsilon)$. \end{definition}

If $\Psi$ is infinitesimally nonnegative, then $\kappa'''(0)\geq 0$ for all commuting $X,Y\in\mg$.  We will use this third derivative restriction to begin classifying the infinitesimally nonnegative endomorphisms $\Psi$.  A subsequent paper will establish that if $\Phi_{t_0}$ has nonnegative curvature for a single $t_0>0$, then $\Psi$ is infinitesimally nonnegative~\cite{SMALL06}, which is a weak version of Conjecture~\ref{ILC}.  Via this subsequent result, all infinitesimal restrictions derived in this paper will imply global restrictions on the family of all left-invariant metrics with nonnegative curvature on $G$. 

Let $a_0<a_1<\cdots<a_l$ be the eigenvalues of $\Psi$, and let $\mp_0,\mp_1,...,\mp_l$ denote the corresponding eigenspaces, which are mutually $h_0$-orthogonal.  The next lemma is our primary tool for deriving rigidity statements; it says that each vector in $\mp_0$ generates other vectors in $\mp_0$.

\begin{lem}\label{k} Assume that $\Psi$ is infinitesimally nonnegative. If $X\in\mp_0$, $Y\in\mg$ and $[X,Y]=0$, then $[X,\Psi Y]\in\mp_0$.
\end{lem}
\begin{proof}
Lemma~\ref{kprime} applied to the span of $X$ and $Y$ gives:
$$(1/6)\kappa'''(0) = a_0|[X,\Psi Y]|^2 - \lb[X,\Psi Y],\Psi [X,\Psi Y]\rb.$$
Since $a_0$ is the smallest eigenvalue, this is negative unless $[X,\Psi Y]\in\mp_0$.
\end{proof}
\section{Known global and infinitesimal metrics on $SO(4)$}
Each known example of a left-invariant metric $h$ with nonnegative curvature on $G=SO(4)$ comes from Cheeger's construction~\cite{Cheeger}.  In this section, we catalog each known example in terms of the eigenvalue and eigenvector structure of the map $\Phi$ representing it with respect to a fixed bi-invariant metric $h_0$, meaning that $h(A,B)=h_0(\Phi A,B)$.  In particular, we will verify that:
\begin{observation} For each known example, $\mg$ decomposes into three 2-dimensional $\Phi$-invariant abelian subalgebras.
\end{observation}
As evidence for Conjecture~\ref{ILC}, we will also verify in this section that:
\begin{observation} For each known left-invariant metric $h$ with nonnegative curvature on $SO(4)$, the unique inverse-linear variation connecting $h_0$ to $h_1=h$ is through nonnegatively curved metrics.
\end{observation}

\vspace{.1in}
\textbf{\emph{(1) PRODUCT METRICS}}.  The Lie algebra $\mg=so(4)$ is a product $\mg=\mg_1\oplus\mg_2$, with each factor isomorphic to $so(3)$.  The two factors are $h_0$-orthogonal.  If they are $h$-orthogonal, then $h$ is a product metric on $SO(4)$'s double cover $S^3\times S^3$.  The classification of product metrics with nonnegative curvature reduces to the classification of left-invariant metrics with nonnegative curvature on $SO(3)$, solved in~\cite{SMALL}.  Observe that for any product metric, $\mg$ decomposes into three 2-dimensional $\Phi$-invariant abelian subalgebras, obtained by pairing eigenvectors from the two factors.

As for infinitesimal examples, we call a self-adjoint map $\Psi:\mg\ra\mg$ a \emph{product} if $\Psi(\mg_1)\subset\mg_1$, which is equivalent to $\Psi(\mg_2)\subset\mg_2$.  Any such $\Psi$ generates a variation $\Phi_t=(I-t\Psi)^{-1}$ through product metrics, which are nonnegative for small $t$, since all product metrics near $h_0$ are nonnegatively curved.  Every product metric $h$ can be reached by such a variation.

It remains to verify that for any left-invariant product metric $h$ with nonnegative curvature, the intermediate metrics of the inverse-linear variation connecting $h_0$ to $h$ are nonnegatively curved.  Clearly it suffices to prove the analogous fact for metrics on $S^3$.

For this, let $g_0$ be a fixed bi-invariant metric on $S^3$.  Let $g_R$ be a right-invariant metric with nonnegative curvature on $S^3$ with eigenvalues $\{\lambda_1,\lambda_2,\lambda_3\}$ and let $a>0$.  Consider the following left-invariant metric on $S^3$:
\begin{equation}\label{above}(S^3,g) = ((S^3,a\cdot g_0)\times(S^3,g_R))/S^3,\end{equation}
where $S^3$ acts diagonally on the right.  The eigenvalues of $g$ are $\left\{\frac{a \lambda_i}{(1+\lambda_i)}\mid i=1,2,3\right\}$.  It was shown in~\cite{SMALL} that this metric is positively curved, and that every positively-curved left-invariant metric on $S^3$ has this form.

Let $\alpha\in\R$.  Replacing $a$ with $1/(1-\alpha t)$ and replacing $g_R$ with $(1/t)g_R$ yields a family $g_t$ of left-invariant metrics on $S^3$.  The eigevalues of $g_t$ are $\left\{\frac{\lambda_i}{(t+\lambda_i)(1-\alpha t)}\mid i=1,2,3\right\}$.  The time $t=0$ derivative of this variation is:
$$\Psi=\left(\begin{matrix} \alpha-1/\lambda_1 & 0 & 0 \\ 0 & \alpha-1/\lambda_2 & 0 \\ 0 & 0 & \alpha-1/\lambda_3\end{matrix}\right).$$
Let $\hat{g}_t$ denote the inverse-linear variation determined by $\Psi$, which is described in matrix form as $\hat{\Phi}_t=(I-t\Psi)^{-1}$.  The eigenvalues of $\hat{g}_t$ are
$\left\{\frac{\lambda_i}{t+\lambda_i(1-\alpha t)}\mid i=1,2,3\right\}$.
If $\alpha=0$, then $\hat{g}_t=g_t$.  If $\alpha\neq 0$, then $g_t$ is not inverse-linear.  Notice that $\hat{g}_t$ can be obtained in equation~\ref{above} by replacing $a$ with $1/(1-\alpha t)$ and $g_R$ with $\left(\frac{1-\alpha t}{t}\right)\cdot g_R$, as long as $1-\alpha t> 0$.  If $\alpha<0$, then $\hat{g}_t$ is nonnegatively curved for all $t$.  If $\alpha>0$, then $\hat{g}_t$ is at least nonnegatively curved until $1-\alpha t=0$, at which time $\hat{g}_t$ has eigenvalues $\{\alpha\lambda_1,\alpha\lambda_2,\alpha\lambda_3\}$.  By choosing $\alpha=1$, we have at this cut-off time $t=1$ that $\hat{g}_1=g_R$.  In summary, $g_R$ was an arbitrary positive curvature metric, and $\hat{g}_t$ is an inverse-linear variation through nonnegatively curved metrics connecting $g_0$ to $g_R$.
\vspace{.1in}

\textbf{\emph{(2) TORUS ACTIONS}}.
A ``Berger metric'' on $\mg_1$ is obtained from $h_0|_{\mg_1}$ by scaling the span of a unit-length vector $A\in\mg_1$ and its orthogonal compliment in $\mg_1$ by separate factors $a$ and $c$ respectively.  Similarly, a unit-length vector $B\in\mg_2$ and scaling factors $b,d>0$ determine a Berger metric on $\mg_2$.  Let $\tilde{h}$ denote the product of these two Berger metrics.  In other words, $\tilde{h}(A,A)=a$, $\tilde{h}(B,B)=b$, and $\tilde{h}=c\cdot h_0$ on $\mg_1\ominus\text{span}\{A\}$, and $\tilde{h}=d\cdot h_0$ on $\mg_2\ominus\text{span}\{B\}$ (where ``$\ominus$'' means the orthogonal compliment).  Considering $\tilde{h}$ as a left-invariant metric on $SO(4)$'s double cover, notice that $(S^3\times S^3,\tilde{h})$ admits an isometric torus action corresponding to the left-invariant Killing fields determined by $A$ and $B$.  We use Cheeger's method to alter the metric via this action:
$$(S^3\times S^3,h)=((S^3\times S^3,\tilde{h})\times (T^2,g_T))/T^2,$$
where $g_T$ is any right-invariant (and hence bi-invariant and flat) metric on the torus $T^2$. 
 
The metric $\tilde{h}$ and hence also $h$ has nonnegative curvature as long as $a/c,b/d\leq 4/3$ (the largest factor by which Hopf circles on $S^3$ can be expanded without losing nonnegative curvature).  

Let $\{e_1,e_2\}$ denote the eigenvectors of $g_T$ with respect to $\tilde{h}$, and let $\{\lambda_1,\lambda_2\}$ denote the eigenvalues, which are arbitrary.  Denote $\tau=\text{span}\{A,B\}$.  The eigenvectors of $h|_{\tau}$ with respect to $\tilde{h}|_{\tau}$ are $\{e_1,e_2\}$.  The eigenvalues are $\{\lambda_1/(1+\lambda_1),\lambda_2/(1+\lambda_2)\}$, which are arbitrary in $(0,1)$.  Thus, $g_T$ can be chosen so as to make $h|_{\tau}$ become any metric which is $\leq\tilde{h}|_\tau$ (which means that $h(X,X)\leq\tilde{h}(X,X)$ for all $X\in\tau$).

In summary, the linear map $\Phi$ describing $h$ in terms of $h_0$ has:
\begin{itemize}
\item One arbitrary eigenvalue $c$ corresponding to a two dimensional subspace of $\mg_1$.
\item One arbitrary eigenvalue $d$ corresponding to a two dimensional subspace of $\mg_2$.
\item Arbitrary eigenvalues and eigenvectors on the remaining subspace $\tau$, subject only to the restriction that the metric on $\tau$ is bounded above by the product metric with scaling factors $(4/3)c$ and $(4/3)d$.
\end{itemize}

Observe that $\mg$ decomposes into three 2-dimensional $\Phi$-invariant abelian subalgebras: one equals $\tau$, and the other two are obtained by pairing vectors in $\mg_1$ with vectors in $\mg_2$.


As for infinitesimal examples, consider a matrix $\Psi$ which has the form
\begin{equation}\label{toruss}\Psi=\left(\begin{matrix}
c & 0 & 0 & 0 & 0 & 0 \\
0 & c & 0 & 0 & 0 & 0 \\
0 & 0 & a_1 & a_3 & 0 & 0\\
0 & 0 & a_3 & a_2 & 0 & 0\\
0 & 0 & 0 & 0 & d & 0 \\
0 & 0 & 0 & 0 & 0 & d
\end{matrix}\right)\end{equation}
with respect to a basis $\{A_1,A_2,A_3,B_1,B_2,B_3\}$, where the $A$'s and $B$'s form orthonormal bases of $\mg_1$ and $\mg_2$ respectively.
Then $\Phi_t=(I-t\Psi)^{-1}$ is a infinitesimally nonnegative variation through the type of ``torus-metrics'' described above.  Each torus-metric can be reached by such a variation.  A calculation shows that if $h$ satisfies the $4/3$ restriction, then so do the intermediate metrics between $h_0$ and $h$.   The metrics of these variations have no more positive curvature planes than $h_0$.  More precisely, for any commuting pair $X,Y\in\mg$, a calculation shows that $\kappa(t)=0$ for all $t$ (with $\kappa$ as in Lemma~\ref{kprime}).

\vspace{.1in}

\textbf{\emph{(3) $S^3$-ACTIONS}}.  Let $\tilde{h}$ denote the bi-invariant metric on $S^3\times S^3$ obtained from $h_0$ by rescaling $\mg_1$ and $\mg_2$ by factors $a$ and $b$ respectively.  Let $g_R$ denote a right-invariant metric with nonnegative curvature on $S^3$ with eigenvalues $\{\lambda_1,\lambda_2,\lambda_3\}$ and eigenvectors $\{e_1,e_2,e_3\}$.  Define:
$$(S^3\times S^3,h)=((S^3\times S^3,\tilde{h})\times(S^3,g_R))/S^3,$$
where $S^3$ acts diagonally.  Consider the basis: $$\mg=\mg_1\oplus\mg_2=\text{span}\{A_1,A_2,A_3\}\oplus\text{span}\{B_1,B_2,B_3\},$$
where $A_i=(e_i,0)$ and $B_i=(0,e_i)$.  Let $V_i=\text{span}\{A_i,B_i\}$, which for each $i$ is a 2-dimensional abelian subalgebra of $\mg$.  Notice that the three $V_i$'s are mutually orthogonal with respect to $h_0$, $\tilde{h}$, and $h$.  It therefore suffices to describe $h$ in terms of $h_0$ separately on each $V_i$.  

For this, the matrix representing $\tilde{h}$ in terms of $h_0$ on $V_i$ in the basis $\{A_i,B_i\}$ is $M_i=\left(\begin{matrix} a & 0 \\ 0 & b\end{matrix}\right)$.  The matrix representing $h$ in terms of $\tilde{h}$ in the basis $\{A_i+B_i,b A_i - a B_i\}$ is $N_i=\left(\begin{matrix} t_i & 0 \\ 0 & 1\end{matrix}\right)$, where $t_i=\frac{\lambda_i}{1+\lambda_i}$. Thus, letting $T$ be the change of basis matrix, $T=\left(\begin{matrix} 1 & b \\ 1 & -a\end{matrix}\right)$, the matrix we seek which represents $h$ in terms of $h_0$ on $V_i$ in the basis $\{A_i,B_i\}$ is:
\begin{equation}\label{willow}\Phi_i= M_i(TN_iT^{-1})=\frac{1}{a+b}\left(\begin{matrix} a(b+at_i) & ab(t_i-1) \\ ab(t_i-1) & b(a+bt_i)\end{matrix}\right).\end{equation}

In summary, $\mg$ decomposes into the three $\Phi$-invariant 2-dimensional abelian subalgebras, $\{V_1,V_2,V_3\}$.  However, with only the five parameters $\{a,b,t_1,t_2,t_3\}$ under our control, and with restrictions on the $t$'s, we do not attain the full 9-parameter family of metrics for which the subalgebras $\{V_1,V_2,V_3\}$ are $\Phi$-invariant.

To produce infinitesimal examples, we consider natural variations through metrics of the type described above.  Choose $\alpha,\beta\in\R$.  In the above construction, replacing $a$ with $1/(1-\alpha t)$, $b$ with $1/(1-\beta t)$, and $g_R$ with $(1/t)g_R$ yields a family $\{h_t\}$ of left-invariant metrics on $S^3\times S^3$ which extends smoothly at $t=0$ to the initial metric $h_0$.  In Equation~\ref{willow}, each $\lambda_i$ gets replaced with $\lambda_i/t$, so each $t_i$ gets replaced with $\lambda_i/(t+\lambda_i)$, and the derivative at $t=0$ of this variation is $\Psi=\text{diag}(\Psi_1,\Psi_2,\Psi_3)$, where:
\begin{equation}\label{gg}\Psi_i=\left(\begin{matrix}\alpha & 0 \\ 0 & \beta\end{matrix}\right) -
\frac{1}{2\lambda_i}\left(\begin{matrix} 1 & 1 \\ 1 & 1\end{matrix}\right).\end{equation}
Let $\{\hat{h}_t\}$ denote the inverse-linear variation determined by $\Psi$, which is described in matrix form as $\hat{\Phi}_t=(I-t\Psi)^{-1}$.
If $\alpha=\beta=0$, then $h_t=\hat{h}_t$, so $\hat{h}_t$ is nonnegatively curved for all $t\geq 0$.  For non-zero $\alpha,\beta$, $\{h_t\}$ is not inverse-linear, but we will nevertheless prove that $\hat{h}_t$ is nonnegatively curved for small $t$.  To establish this, we must find parameters $\{\overline{a},\overline{b},\overline{\lambda}_1\overline{\lambda}_2,\overline{\lambda}_3\}$ such that $$\hat{\Phi}_t^i=\frac{1}{\overline{a}+\overline{b}}\left(\begin{matrix} \overline{a}(\overline{b}+\overline{a}\overline{t}_i) & \overline{a}\overline{b}(\overline{t}_i-1) \\ \overline{a}\overline{b}(\overline{t}_i-1) & \overline{b}(\overline{a}+\overline{b}\overline{t}_i)\end{matrix}\right)$$
for each $i=1,2,3$, where $\hat{\Phi}_t^i$ denotes the $i^{\text{th}}$ diagonal block of $\hat{\Phi}_t$, and $\overline{t}_i=\overline{\lambda}_i/(1+\overline{\lambda}_i)$.  A solution to this system of equations is:
\begin{equation}\label{yy}\overline{a}=\frac{1}{1-\alpha t},\,\,
  \overline{b}=\frac{1}{1-\beta t},\,\,
  \overline{\lambda}_i=
       2\lambda_i\frac{(1-\alpha t)(1-\beta t)}{(1-\alpha t)+(1-\beta t)},\end{equation}
as long as $1-\alpha t$ and $1-\beta t$ are positive.  
Since $(\overline{\lambda}_1,\overline{\lambda}_2,\overline{\lambda}_3)$ is a positive multiple of $(\lambda_1,\lambda_2,\lambda_3)$, $\hat{h}_t$ is nonnegatively curved as long as $1-\alpha t$ and $1-\beta t$ are positive.  It is easy to see from equation~\ref{yy} that every known metric of the above form can be achieved as $\hat{h}_1$ by proper choice of $\{\alpha,\beta,\lambda_1,\lambda_2,\lambda_3\}$ with $\alpha,\beta<1$.
\section{Infinitesimal rigidity for $SO(4)$}
In this section, we assume that $G=SO(4)$ and $\Psi:\mg\ra\mg$ is infinitesimally nonnegative, and we prove rigidity results for $\Psi$, using only the property that $\kappa'''(0)\geq 0$ for all commuting $X,Y\in\mg$.  

Recall that $\mg = so(4)=\mg_1\oplus\mg_2$ is a product, and $X\in\mg$ is called \underline{regular} if it has non-zero projections onto both $\mg_1$ and $\mg_2$; otherwise, it is called \underline{singular}.

The previous section classified the known possibilities of $\Psi$ into three types, coming from: (1) products, (2) torus actions and (3) $S^3$-actions.  In the first two cases, $\Psi$ has a non-zero singular eigenvector, while in the third case, it does not.

\begin{theorem}\label{Th1} If $\Psi$ has a non-zero singular eigenvector, the either $\Psi$ is a product or $\Psi$ has the form of Equation~\ref{toruss}.  In either case, $h_t$ is a family of known examples with nonnegative curvature for sufficiently small $t$.\end{theorem}

If $\Psi$ has no non-zero singular eigenvectors, we conjecture that $\Psi$ is a known example coming from an $S^3$-action.  A first step in this direction is to locate three $\Psi$-invariant abelian subalgebras.  The following theorem falls just short of this goal:
\begin{theorem}\label{Th2}
There are orthonormal bases $\{A_1,A_2,A_3\}$ and $\{B_1,B_2,B_3\}$ of the two factors of $\mg=\mg_1\oplus\mg_2$ such that with respect to the basis $\{A_1,B_1,A_2,B_2,A_3,B_3\}$, $\Psi$ has the form:
$$\Psi=\left(\begin{matrix}
a_1 & a_3 & 0 & 0 & 0 & 0 \\
a_3 & a_2 & 0 & 0 & 0 & 0 \\
0 & 0 & b_1 & b_3 & \lambda & 0\\
0 & 0 & b_3 & b_2 & 0 & \mu\\
0 & 0 & \lambda & 0 & c_1 & c_3 \\
0 & 0 & 0 & \mu & c_3 & c_2
\end{matrix}\right).$$
\end{theorem}
We conjecture that $\lambda=\mu=0$, which means that $\mg$ decomposes into three orthogonal $\Psi$-invariant abelian subalgebras, as it should.  Even granting this conjecture, there remains the work of reducing the above 9-parameter family to the 5-parameter family of known examples from Equation~\ref{gg}.  This appears to be a computationally difficult problem.

The remainder of this paper is devoted to proving Theorems~\ref{Th1} and~\ref{Th2}.  We begin with a weak version of Theorem~\ref{Th1}.  Recall that $\mp_0$ denotes the eigenspace corresponding to the smallest eigenvalue, $a_0$, of $\Psi$. 
\begin{lem}\label{nongen}If $\mp_0$ contains a non-zero singular vector, then either $\Psi$ is a product or $\Psi$ has the form of Equation~\ref{toruss}.
\end{lem}
\begin{proof}
Without loss of generality, assume there exists a non-zero vector $X_1\in\mg_1\cap\mp_0$.  Assume that $\Psi$ is not a product, so there exists $\hat{Y}\in\mg_2$ such that $\Psi\hat{Y}$ has a nonzero projection, $X_2$, onto $\mg_1$.  Notice that $X_1$ and $X_2$ are orthogonal because:
$$\lb X_1,X_2\rb = \lb X_1,\Psi \hat{Y}\rb = \lb \Psi X_1,\hat{Y}\rb = a_0\lb X_1,\hat{Y}\rb = 0.$$
Let $X_3=[X_1,\Psi\hat{Y}]\in\mg_1$, which by Lemma~\ref{k} lies in $\mp_0$, so $\text{span}\{X_1,X_3\}\subset\mp_0$.  Let $Y_2=\Psi (X_2)$, which is a non-zero vector in $\mg_2$ (by the self-adjoint property of $\Psi $). Complete $\{Y_2\}$ to an orthogonal basis $\{Y_1,Y_2,Y_3\}$ of $\mg_2$.  Notice that $\Psi (\text{span}\{Y_1,Y_3\})\subset\mg_2$ (again by the self-adjoint property of $\Psi $).  In summary, after scaling all the vectors to unit-length, we have an orthonormal bases:
$$\mg = \mg_1\oplus\mg_2 = \text{span}\{X_1,X_2,X_3\}\oplus\text{span}\{Y_1,Y_2,Y_3\}$$
with $\text{span}\{X_1,X_3\}\subset\mp_0$, and $\Psi  X_2 = cY_2+\lambda X_2$ (for some $c,\lambda\in\R$ with $c\neq 0$), and $\Psi (\text{span}\{Y_1,Y_3\})\subset\mg_2$.

Applying Lemma~\ref{kprime} to the vectors $X_2$ and $Y_1$ gives:
\begin{eqnarray*}
\kappa'''(0) & = & 6\lb[\Psi  X_2,Y_1],[\Psi X_2,\Psi Y_1]\rb
             - 6\lb [\Psi  X_2,Y_1],\Psi [\Psi  X_2,Y_1]\rb \\
        & = & 6\lb [ cY_2,Y_1],[cY_2,\Psi  Y_1]\rb - 6\lb[cY_2,Y_1],\Psi [cY_2,Y_1]\rb\\
        & = & -6c^2\lb Y_3,[Y_2,\Psi  Y_1]\rb - 6c^2\lb Y_3,\Psi  Y_3\rb \geq 0.
\end{eqnarray*}
Notice that:
\begin{eqnarray*}
\lb Y_3,[Y_2,\Psi  Y_1]\rb & = & \lb Y_3,[Y_2,\text{projection of }\Psi Y_1\text{ onto }Y_1]\rb \\
                             & = & \lb Y_3,[Y_2,\lb\Psi Y_1,Y_1\rb Y_1]\rb \\
                             & = & -\lb \Psi Y_1,Y_1\rb.
\end{eqnarray*}
The conclusion is:
$$\lb Y_1,\Psi Y_1\rb\geq\lb Y_3,\Psi Y_3\rb.$$
We're assuming above that $\{Y_1,Y_2,Y_3\}$ is oriented such that $[Y_1,Y_2]=Y_3$, etc., but the opposite orientation assumption yields the same conclusion.

Similarly, applying Lemma~\ref{kprime} to the vectors $X_2$ and $Y_3$ yields the reverse inequality, so:
$$\lb Y_1,\Psi Y_1\rb = \lb Y_3,\Psi Y_3\rb.$$

Replacing $Y_1$ and $Y_3$ with any other orthonormal basis of $\text{span}\{Y_1,Y_3\}$ yields the same conclusion.  In other words, for any angle $\theta$, set $a=\cos(\theta)$ and $b=\sin(\theta)$, and:

$$\lb aY_1+bY_3,\Psi (aY_1+bY_3)\rb = \lb bY_1-aY_3,\Psi (bY_1-aY_3)\rb.$$
This implies that $\lb Y_1,\Psi Y_3\rb = \lb \Psi Y_1,Y_3\rb = 0$.  The linear map from $\text{span}\{Y_1,Y_3\}$ to $\R$ sending $Y\mapsto\lb \Psi Y,Y_2\rb$ has a non-zero vector in its kernel.  Assume without loss of generality that $Y_1$ is in its kernel.  Notice that $Y_1$ is an eigenvector of $\Psi $.

In the ordered basis $\{X_1,X_2,X_3,Y_1,Y_2,Y_3\}$, we thus far have:
$$\Psi =\left(\begin{matrix}
a_0 & 0 & 0 & 0 & 0 & 0 \\
0 & \lambda & 0 & 0 & c & 0 \\
0 & 0 & a_0 & 0 & 0 & 0\\
0 & 0 & 0 & \beta & 0 & 0\\
0 & c & 0 & 0 & \gamma & s \\
0 & 0 & 0 & 0 & s & \beta
\end{matrix}\right)$$
Applying our $\kappa'''(0)$ formula to $X=X_2$ and $Y=aY_2+bY_3$ gives:
$$\kappa'''(0) = 6bc^2(as+b\beta) - 6 b^2c^2\beta = 6bc^2as.$$
Since $\kappa'''(0)\geq 0$ for all choices of $\{a,b\}$, and $c\neq 0$, we learn that $s=0$.  After re-ordering the basis, $\Psi$ has the form of Equation~\ref{toruss}.
\end{proof}
\begin{theorem}
The eigenspace $\mp_0$ contains a non-zero vector which belongs to a $\Psi$-invariant $2$-dimensional abelian subalgebra of $\mg$.
\end{theorem}
\begin{proof}
If $\mp_0$ contains a non-zero singular vector, the conclusion follows easily from Lemma~\ref{nongen}, so we assume that this is not the case.  When $A=(A_1,A_2)\in\mg=\mg_1\oplus\mg_2$ is regular, let $\overline{A}=\left(\frac{|A_2|}{|A_1|}A_1,-\frac{|A_1|}{|A_2|}A_2\right)$, which commutes with $A$, is orthogonal to $A$, and has the same norm as $A$.

The proof is indirect.  We assume for each $A\in\mp_0$ that $\text{span}\{A,\overline{A}\}$ is not $\Psi$-invariant, and we derive a contradiction.

Let $A\in\mp_0$ be unit-length.  Since $\Psi $ is self-adjoint, $\Psi \overline{A}$ is orthogonal to $A$.  Notice that $\overline{A}$ is not an eigenvector of $\Psi $; if it were, then $\text{span}\{A,\overline{A}\}$ would be an invariant abelian subalgebra.  Therefore, $[A,\Psi \overline{A}]$ is non-zero.  Let $B$ be the unit-length vector in the direction of $[A,\Psi \overline{A}]$.  By Lemma~\ref{k}, $B\in\mp_0$.  Notice that $B$ is orthogonal to $A$ and $\overline{A}$.

So far we know that $\dim(\mp_0)\geq 2$.  Clearly $\dim(\mp_0)\leq 3$ because it contains no non-zero singular vectors, and hence intersects $\mg_1$ and $\mg_2$ trivially.  We wish to prove $\dim(\mp_0)=2$.  Suppose to the contrary that $\dim(\mp_0)=3$.  Consider the map from $\mp_0$ to $\mp_0$ defined as
$$Z\mapsto [Z,\Psi\overline{Z}].$$
By the above arguments, this map sends each unit-length $Z\in\mp_0$ to a non-zero vector in $\mp_0$ orthogonal to $Z$.  This map therefore induces a smooth non-vanishing vector field on the unit $2$-sphere in $\mp_0$, which is a contradiction.  Thus, $\dim(\mp_0)=2$.  Notice $A$ and $B$ play symmetric roles in that $[B,\Psi\overline{B}]$ is parallel to $A$ (because it lies in $\mp_0$ and is perpendicular to $B$), and $A$ is orthogonal to $B$ and $\overline{B}$.

Choose unit-length vectors $C_1\in\mg_1$ and $C_2\in\mg_2$ such that $\{A,\overline{A},B,\overline{B},C_1,C_2\}$ is an orthonormal basis of $\mg$.  For $i=1,2$, the $\mg_i$-components of $\{A,B,C_i\}$ form an orthogonal basis of $\mg_i$.  The $C_i$'s can be chosen so that these orthogonal bases are oriented, so after normalizing, they act like $\{\ii,\jj,\kk\}$ with respect to their Lie bracket structure.  For purposes of calculating Lie brackets in this basis, we lose no generality in assuming that for some $a,b\in(0,1)$,
\begin{gather}\label{iii}
A=(a\ii,\sqrt{1-a^2}\ii), \,\,\,\, B=(b\jj,\sqrt{1-b^2}\jj),\,\,\,\, C_1=(\kk,0) \\
\overline{A}=(\sqrt{1-a^2}\ii,-a\ii), \,\,\,\, \overline{B}=(\sqrt{1-b^2}\jj,-b\jj),\,\,\,\, C_2=(0,\kk).\notag
\end{gather}
Notice that $h_0(\Psi\overline{A},\overline{B})=h_0(\Psi\overline{B},\overline{A})=0$, because if $\Psi\overline{A}$ had a nonzero $\overline{B}$-component, then $[A,\Psi\overline{A}]$ would have nonzero $C_1$ and $C_2$-components.

In the basis $\{A,\overline{A},B,\overline{B},C_1,C_2\}$, $\Psi $ has the form:
\begin{equation}\label{early}\Psi =\left(\begin{matrix}
a_0 & 0 & 0 & 0 & 0 & 0 \\
0 & p & 0 & 0 & \alpha_1 & \alpha_2 \\
0 & 0 & a_0 & 0 & 0 & 0\\
0 & 0 & 0 & q & \beta_1 & \beta_2\\
0 & \alpha_1 & 0 & \beta_1 & f_1 & f_2 \\
0 & \alpha_2 & 0 & \beta_2 & f_2 & f_3
\end{matrix}\right).\end{equation}
There are a few obvious restrictions among the variables determining $\Psi $.  For example, since $[A,\Psi \overline{A}]$ is parallel to $B$, and $[B,\Psi \overline{B}]$ is parallel to $A$, we learn:
\begin{equation}\label{simon}\frac{\alpha_1}{\alpha_2}=\frac{\beta_2}{\beta_1}
=\frac{b\sqrt{1-a^2}}{a\sqrt{1-b^2}}.\end{equation} and we obtain:

\begin{equation}\label{early2}\Psi =\left(\begin{matrix}
a_0 & 0 & 0 & 0 & 0 & 0 \\
0 & p & 0 & 0 & \alpha & \alpha \cdot s \\
0 & 0 & a_0 & 0 & 0 & 0\\
0 & 0 & 0 & q & \beta \cdot s & \beta\\
0 & \alpha & 0 & \beta \cdot s & f_1 & f_2 \\
0 & \alpha \cdot s & 0 & \beta & f_2 & f_3
\end{matrix}\right),\end{equation}
where $s = \frac{ a \sqrt{1-b^{2}} }{ b \sqrt{1-a^{2}} } > 0$ and $\alpha, \beta \neq 0$.

Using Lemma~\ref{k}, we can now prove that $s=1$ and consequently $a = b$. Indeed,  for every $Z\in\text{span}\{A,B\}$, we have $[Z,\Psi \overline{Z}]\in\text{span}\{A,B\}$.  In particular, let $Z_t=(\cos t) A+(\sin t)B$, so:
$$\overline{Z}_t = \left(  f(t) \left(a\cos(t) \ii + b\sin(t) \jj\right),  - (1/f(t)) \left( \sqrt{1-a^{2}} \cos(t) \ii + \sqrt{ 1-b^{2}}\sin(t)\jj \right)\right),$$
where 
$$f(t)=  \sqrt{\frac{ (1-a^{2})\cos^{2}(t) + (1-b^{2})\sin^{2}(t) }{ a^{2}\cos^{2}(t) + b^{2}\sin^{2}(t)} }.$$
We will use that the following vector lies in $\text{span}\{A,B\}$:
\begin{eqnarray*}
Q &=& \ddt\Big|_{t=0}[Z_t,\Psi\overline{Z_t}] = [B,\Psi\overline{A}] + \left[A, \Psi \left(f'(0)a\ii + f(0) b \jj, -g'(0) \sqrt{1-a^{2}}  \ii -g(0) \sqrt{ 1-b^{2}} \jj\right)\right] \\
    &=& [B,\Psi\overline{A}] + \left[A, \Psi \left(f(0) b \jj, -g(0) \sqrt{ 1-b^{2}} \jj\right)\right] \\
    &=& [B,\Psi\overline{A}] + \left[A, \Psi \left( \frac{b \sqrt{1-a^{2}}}{a} \jj , -\frac{a \sqrt{1-b^{2}}}{\sqrt{ 1-a^{2} }} \jj \right)\right] \\
    &=& [B,\Psi\overline{A}] + \left[A, \Psi \left( \sqrt{1-b^{2}} \cdot  s^{-1} \, \jj , -b \cdot  s \, \jj\right)\right].
\end{eqnarray*}
In particular, $Q$ is perpendicular to $\overline{A}$, so:
\begin{eqnarray*}
0 &=& \lb Q,\overline{A}\rb = \lb [B,\Psi\overline{A}],\overline{A}\rb + 
\left\langle\left[A, \Psi \left( \sqrt{1-b^{2}} \cdot  s^{-1} \, \jj , -b \cdot  s \, \jj\right)\right],\overline{A}\right\rangle \\
   &=&  \lb [B,\Psi\overline{A}],\overline{A}\rb = - \langle \Psi \overline{A}, [B, \overline{A} ] \rangle \\
   &=& - \langle p \overline{A} + (\alpha\kk , \alpha s\kk) , [(b\jj, \sqrt{1-b^{2}} \jj), (\sqrt{1-a^{2}} \ii, -a\ii)] \rangle \\
   &=& - \langle p \overline{A} + (\alpha\kk , \alpha s\kk) , ( - b \sqrt{1-a^{2}} \kk, a \sqrt{1-b^{2}} \kk) \rangle \\
   &=& \alpha b \sqrt{1-a^{2}} - s\alpha a \sqrt{1-b^{2}},\\
\end{eqnarray*}
which implies that $s = \frac{b \sqrt{1-a^{2}}}{a \sqrt{1-b^{2}} } = s^{-1}$. It follows that $s=1$ and, consequently, $a= b$. Now the fact that the orthogonal projection of $Q$ onto $\text{span}\{C_1,C_2\}$ is zero is equivalent to 
\begin{equation}\label{gab}p(-b\sqrt{1-a^2}\kk,a\sqrt{1-b^2}\kk)+q(a\sqrt{1-b^2}\kk,-b\sqrt{1-a^2}\kk)=0.\end{equation} 
Since $a=b$, this implies that $q=p$. So we obtain 
\begin{equation}\label{psi}\Psi =\left(\begin{matrix}
a_0 & 0 & 0 & 0 & 0 & 0 \\
0 & p & 0 & 0 & \alpha & \alpha \\
0 & 0 & a_0 & 0 & 0 & 0\\
0 & 0 & 0 & p & \beta & \beta\\
0 & \alpha & 0 & \beta & f_1 & f_2 \\
0 & \alpha & 0 & \beta & f_2 & f_3
\end{matrix}\right).\end{equation}
Since $a=b$, it is easy to see that $[A,\overline{B}]+[B,\overline{A}]=0$.
This implies that $V_1=\beta\overline{A}-\alpha\overline{B}$ commutes with $V_2=\beta A-\alpha B.$  Since $V_2\in\mp_0$, and $V_1$ is an eigenvector of $\Psi$ (with eigenvalue $p$), we learn that $\text{span}\{V_1,V_2\}$ is a $\Psi$-invariant 2-dimensional abelian subalgebra of $\mg$ containing a non-zero vector in $\mp_0$.  This is a contradiction.
\end{proof}

\begin{proof}[Proof of Theorem~\ref{Th2}]
By the previous theorem, there exists a $\Psi$-invariant abelian subalgebra of $\mg$, spanned by some $A_1\in\mg_1$ and some $B_1\in\mg_2$.  Let $V_1$ denote the orthogonal compliment of $A_1$ in $\mg_1$, and let $V_2$ denote the orthogonal compliment of $B_1$ in $\mg_2$.

Let $\pi_1:\mg\ra\mg_1$  and $\pi_2:\mg\ra\mg_2$ denote the projections.  Define $T_1:V_1\ra V_2$ as $T_1=\pi_2\circ\Psi|_{V_1}$, and define $T_2:V_2\ra V_1$ as $T_2=\pi_1\circ\Psi|_{V_2}$.  Notice that for all $A\in V_1$ and $B\in V_2$,
$$\lb T_1A,B\rb = \lb\Psi A,B\rb = \lb A,\Psi B\rb = \lb A,T_2 B\rb.$$

Let $S^1$ denote the circle of unit-length vectors in $V_1$.  Let $R:S^1\ra S^1$ denote a $90^\circ$ rotation.  Define $F:S^1\ra\R$ as $F(A)=\lb T_1(A),T_1(R(A))\rb$.  Notice that for all $A\in S^1$,
$$F(R(A))=\lb T_1(R(A)),T_1(-A)\rb = -F(A).$$
This implies that there exists $A_2\in S^1$ such that $F(A_2)=0$.  Let $A_3=R(A_2)$.  First suppose that $T_1$ (and hence also $T_2$) is nonsingular.  Define $B_2=T_1(A_2)/|T_1(A_2)|$ and $B_3=T_1(A_3)/|T_1(A_3)|$.  The fact that $F(A_2)=0$ immediately implies that $B_2$ and $B_3$ are orthogonal, and that $T_2(B_2)\parallel A_2$ and $T_2(B_3)\parallel A_3$.  Thus, the basis $\{A_1,A_2,A_3,B_1,B_2,B_3\}$ satisfies the conclusion of the theorem.

If $T_1$ (and hence also $T_2$) is singular, then arbitrary orthonormal bases $\{A_2,A_3\}$ of $V_1$ and $\{B_2,B_3\}$ of $V_2$ work, so long as $A_2\in\text{ker}(T_1)$ and $B_2\in\text{ker}(T_2)$.
\end{proof}
Our final proof is due to Nela Vukmirovic and Zachary Madden:
\begin{proof}[Proof of Theorem~\ref{Th1}]
Choose bases $\{A_1,A_2,A_3\}$ of $\mg_1$ and $\{B_1,B_2,B_3\}$ of $\mg_2$ so that $\Psi$ has the matrix form of Theorem~\ref{Th2}.  If $a_3=0$, then the result follows from Lemma~\ref{nongen}, so we can assume that $a_3\neq 0$.  The hypothesis that $\Psi$ has a non-zero singular eigenvectors implies that $b_3=0$ or $c_3$=0.  Without loss of generality, assume $b_3=0$.

Henceforth, the value $\kappa'''(0)$ with respect to the commuting pair $X=\alpha_1 A_1+\alpha_2 A_2+\alpha_3 A_3$ and $Y=\beta_1 B_1+\beta_2 B_2+\beta_3 B_3$ will be denoted as $[\alpha_1,\alpha_2,\alpha_3,\beta_1,\beta_2,\beta_3]$.

First, $[0,\pm 1,1,1,0,0]=c_3^2(a_2-b_2)\pm 4a_3^2\lambda\geq 0$.  However, since $[0,0,1,0,1,0] + [0,0,1,0,0,1] = c_3^2(b_2-a_2) \geq 0$, we can deduce that $\lambda =0$ and consequently that $c_3^2(b_2-a_2)=0$.

Similarly, $[1,0,0,0,\pm 1,1] = c_3^2(a_1-b_1)\pm 4a_3^2\mu \geq 0$. However, since $[1,0,0,0,1,0] + [1,0,0,0,0,1] = c_3^2(b_1-a_1) \geq 0$ it follows that $\mu =0$ and that $c_3^2(b_1-a_1)=0$.

Further, $[0,1,0,1,0,0]\geq 0$ and $[0,0,1,1,0,0]\geq 0$ give respectively the plus and minus versions of: $ \pm a_3^2(b_1-c_1) \geq 0$.  Similarly from $[1,0,0,0,1,0]$ and  $[1,0,0,0,0,1]$ we conclude that $\pm a_3^2(b_2-c_2) \geq 0$.
Since $a_3$ is non-zero we get that $b_1 = c_1$ and $b_2=c_2$.

Finally, suppose that $c_3$ is not equal to zero. Then from previous equations it follows that $a_1=b_1$ and $a_2=b_2$. By considering $[1,1,1,1,1,1]$, $[1,1,1,-1,1,1]$,
$[1,1,1,1,-1,1]$, and $[1,1,1,1,1,-1]$ we deduce that $\pm a_3^2c_3 \geq 0$, which implies that $c_3=0$.  Thus, $\Psi$ has the form of Equation~\ref{toruss}.
\end{proof}

\bibliographystyle{amsplain}

\end{document}